\documentclass{amsart}
\usepackage{amsmath,amssymb,amsthm,graphicx,epsfig,float,url}
\usepackage[colorlinks=true,citecolor={Plum},linkcolor={Periwinkle}]{hyperref}
\usepackage{pdfsync}
\usepackage[usenames, dvipsnames]{xcolor}
\usepackage{tikz}
\usepackage{mathtools,physics}
\usepackage{subfig}
\usepackage{stmaryrd}
\usepackage{mathrsfs}  
\usepackage{esint}
\usepackage{verbatim}
\usepackage{geometry}
\usepackage{algorithm,algpseudocode} 
\usetikzlibrary{decorations.pathmorphing}
\usetikzlibrary{patterns.meta, plotmarks,calc}

\usepackage{dsfont}
\definecolor{cb-blue}       {RGB}{ 0, 109, 219}
\definecolor{cb-green-lime} {RGB}{ 36, 255,  36}

\newcommand{\R}{\mathbb R}

\def\O{{\Omega}}
\def\n{{\nabla}}
\def\p{{\varphi}}

\def\fooaux#1#2{%
  \mkern2mu
  \setbox0=\hbox{\mathsurround=0pt$#1\bar{\mkern-2mu #2 \mkern2mu}$}%
  \setbox1=\hbox to \wd0{\hss$#1\cdot$\hss}%
  \vbox{\offinterlineskip\copy1\vskip-.4\ht1\box0}%
  \mkern-2mu
}

\newtheorem{theorem}{Theorem}

\newtheorem{material}{material}[section]
\newtheorem{proposition}[material]{Proposition}
\newtheorem{corollary}[material]{Corollary}

\newtheorem{lemma}[material]{Lemma}

\newtheorem{remark}[material]{Remark}

\numberwithin{equation}{section}

\title{Uniqueness \& Attractiveness of stationary solutions in a harvesting Mean Field Game}

\author{Greta Lamonaca}
\address{Greta Lamonaca, Universit\'e d'Orl\'eans, Universit\'e de Tours, CNRS, IDP, UMR 7013, Orl\'eans, France, CEREMADE, UMR CNRS 7534, Universit\'e Paris-Dauphine, Universit\'e PSL, Place du Mar\'echal De Lattre De Tassigny, 75775 Paris cedex 16, France} \email{greta.lamonaca@univ-orleans.fr}
\author{ Idriss Mazari-Fouquer}
\address{Idriss Mazari-Fouquer, CEREMADE, UMR CNRS 7534, Universit\'e Paris-Dauphine, Universit\'e PSL, Place du Mar\'echal De Lattre De Tassigny, 75775 Paris cedex 16, France}
\email{mazari@ceremade.dauphine.fr}
\author{Gr\'egoire Nadin}
\address{Gr\'egoire Nadin, Universit\'e d'Orl\'eans, Universit\'e de Tours, CNRS, IDP, UMR 7013, Orl\'eans, France}
\email{gregoire.nadin@cnrs.fr}

\begin{document}
\begin{abstract}
   Sustainable harvesting of natural resources is a central challenge both in ecology and economics. This paper explores a MFG model for harvesting in an aquatic ecosystem. Such systems have recently attracted attention from the applied mathematics community, and, due to the interplay between the competition of agents and the internal dynamics of the fish population, the analytical properties of the resulting system are very delicate to study. In the present work, we focus on stationary solutions, with a strong emphasis on their uniqueness and dynamical stability.
\end{abstract}

\maketitle


\textbf{AMS-Subject Classification:} 35K10, 35K57, 35Q89, 35Q92, 49N80, 91A12, 91A16.

\textbf{Acknowledgement:}  This project has received funding from the European Union's Horizon Europe research and innovation program under the Marie Sk\l odowska-Curie grant agreement No 101126554.  I. M.-F. is supported by the PSL Young Researcher Starting Grant (P.I: I. M.-F.) ``Optimal control of ecological systems". G. N. and G. L. are partially supported by the ANR project ReaCh ``R\'eaction-diffusion: nouveaux Challenges". I. M.-F. and G. N. are partially supported by the ANR project STOIQUES.
\textbf{Disclaimer:} Co-Funded by the European Union. Views and opinions expressed are however those of the author(s) only and do not necessarily reflect those of the European Union. Neither the European Union nor the granting authority can be held responsible for them. Generative AI (ChatGPT6) was used in the last stages of writing the papers to check for typos or streamline proofs. Authors take full responsibility for the contents of the paper. 

\section{Introduction}
\subsection{Scope of the paper}

The goal of this paper is to analyse models arising in the optimal management of fisheries, where an infinite number of (selfish) agents are competing to maximise their harvest; of particular importance to us are questions related to the long-time behaviour of these scenarios. The underlying paradigm is that of Mean-Field Games (MFG) systems. This work is part of a broader effort to understand rational harvesting and the impact of competition on ecosystems
\cite{Braverman2009,bressan2019competitive,Bressan_2013,zbMATH07248633,zbMATH07006536}, but also fits within two lines of research: obtaining a better understanding of the impact of fisheries on ecosystems \cite{BBC,BBC2,BBCwaste, BBC50} and studying the analytical properties of MFG models\cite{CLLP12,CLLP13,MKB23,zbMATH08109759}.  

While we defer a precise presentation of the model to the next paragraph, let us point out that such MFG models were introduced, and their analysis started, in \cite{zbMATH08109759}. A characteristic of the resulting system is the mix between the standard MFG system (\emph{i.e.} a combination of Hamilton-Jacobi and Fokker-Planck equations describing the value function and evolution of the population of fishermen respectively) and reaction-diffusion equations (describing the evolution of the fish population). A core question in MFG systems is that of \emph{uniqueness} of solutions, and the \emph{convergence to stationary states}, the latter solving the so-called ``ergodic system".  In \cite{zbMATH08109759}, uniqueness and dynamical stability of solutions were proved in perturbative regimes, essentially assuming the interactions between agents to be very weak. In our previous work \cite{arXiv:2606.30106}, we investigated a different situation, that of a coordinated fishery, when a central planner instructs all fishermen as to how to behave and we showed that while, in some perturbative settings, uniqueness of solutions still held, one should expect, depending on biological parameters of the problem, the emergence of multiple solutions. In the present work, moving back to the competitive setting, we investigate whether one can obtain tractable conditions for uniqueness of stationary solutions and their stability, beyond the aforementioned perturbative setting. 

Summarising, the core of this paper is the derivation of sufficient conditions for uniqueness and stability of ergodic solutions. All works regarding such analytical properties rely  either on smallness assumptions (\emph{e.g.} smallness of the time horizon) that are unsuitable for long-time behaviour,  on structural properties (some type of convexity) or on some functional inequalities (dubbed Lasry-Lions monotonicity conditions), the latter of which rely on the comparison between any pair  of solutions. Here, we show that controlling only \emph{one} particular solution can lead to the desired conclusion: uniqueness and dynamical stability of the stationary state. To the best of our knowledge, this is the first such criterion to be derived. We go further by providing several examples where this condition is met. All these constructions rely on fine estimates for the MFG system and the reaction-diffusion model.

\textbf{Model under consideration.} Following \cite{MKB23, zbMATH08109759}, we consider a population of fishermen harvesting a fish population and engaged in  a competition to maximise their payoff. We consider a smooth bounded domain $\O\subset \R^n$, $n\geq 1$ (when $n=1$, we take $\O=(0;\ell)$). $\O$ will be dubbed ``the lake". The fish population is modelled through a density $\theta$ and, in the absence of fishermen, is subject to a Brownian diffusion with strength $\sqrt{2d}$ ($d>0$), to a Malthusian death rate $-\theta^2$ and can access resources $k\in L^\infty(\O)$, leading to a growth rate $k\cdot\theta$. Fix a time horizon $T>0$. For a given  density of fishermen $m\in L^\infty((0;T)\times \O;\R^+)$  and some initial condition $\theta_0\geq 0,\, \not\equiv 0$, the overall evolution of the fish population is governed by the monostable equation 
\[\begin{cases}\partial_t\theta_m-d\Delta\theta_m=\theta_m\left(k-m-\theta_m\right), \ \theta_m >0 &\text{ in }(0;T)\times \O,\,
\\ \partial_\nu \theta_m=0&\text{ on }(0;T)\times \partial \O,\, 
\\ \theta_m(t=0,\cdot)=\theta_0(\cdot)&\text{ in }\O.\end{cases}\]
We assume that the fisherman can control his/her velocity, which is also subject to some random diffusion with strength $\sqrt{2\nu}$, $\nu>0$, so that its trajectory $\left(X_t\right)_{t>0}$ satisfies the following dynamic:
\[ dX^\alpha_t= \alpha_t dt + \sqrt{2\nu}dB_t,\]
where $\left(B_t\right)_{t\ge 0}$ is a standard  Brownian motion, $\alpha$ is the control and the stochastic process is reflected at the boundary of the domain. The strategy $\alpha$ is chosen so as to solve the optimisation problem
\begin{equation}\label{Eq:Payoff}
    \sup_{\alpha \in L^2((0;T)\times \O)} \mathbb{E}\left[\int_0^T \theta_{m} (t, X^\alpha_t)  -\frac{\abs{\alpha(t, X^\alpha_t)}^2}{2}  \ dt \right],
\end{equation}
 so that a single player optimises his/her path while the rest of the population follows a pre-determined evolution. An MFG equilibrium corresponds to a fixed point solution of this optimisation, in the sense that the population density generated by the optimal responses coincides with the prescribed density.
  \begin{remark}Clearly, the underlying assumption that each player has a negligible impact on the fish population may be unrealistic in the context of fishing, since the technology now available to a single fisherman could significantly affect the persistence or extinction of the fish population. Nevertheless, it provides a simplified framework that allows for a tractable analysis of the competitive scenario. \end{remark} Introducing the value function \[ u(t,x)= \sup_{\alpha \in L^2((0;T)\times\O)} \mathbb{E}\left[\int_t^T \theta_{m}(s, X^\alpha_s)   -\frac{\abs{\alpha(s,X_s^\alpha)}^2}{2}  \ ds \right],\] the Dynamic Programming Principle shows, at least formally, that $u$ solves a backward Hamilton-Jacobi equation, and that $\n u$ provides an optimal control for the problem. A MFG solution, corresponding to fixed points in the optimisation problem above, can thus be characterised through the following MFG system (where each equation is endowed with Neumann boundary conditions):\begin{equation}\label{MFG}
    \begin{cases}
        -\partial_t u -\nu \Delta u - \frac{\abs{\nabla u}^2}{2}=\theta \quad & \text{in } (0,T)\times\Omega,\\
        \partial_t m -\nu \Delta m +\n \cdot (m \, \nabla u)=0, \ \fint_\Omega m =1, \ m>0 \ &\text{in } (0,T)\times\Omega,\\
        \partial_t \theta -d \Delta \theta -\theta (k(x)-m-\theta)=0, \ \theta>0 &\text{in } (0,T)\times\Omega,\\
        u(T, \cdot)=0, \ m(0,\cdot)=m_0(\cdot)>0, \ \fint_\Omega m_0=1 &\text{in } \Omega,\\
        \theta(0,\cdot)=\theta_0(\cdot)\not\equiv 0, \ge 0.
    \end{cases}
\end{equation}  
We do not give in this article a rigorous stochastic control justification of the correspondence between optimal controls and solutions to the MFG system.

The basic questions concerning \eqref{MFG} are existence and uniqueness of solutions, and existence, uniqueness and stability of stationary solutions. The existence of solutions to the evolution equation is  guaranteed by a standard fixed point argument which we do not detail here, under the assumption that $\fint_\O k>1$ (which will be one of our working hypotheses). Our main question regards the uniqueness and stability of such steady states. Let us consider the rescaled functions: 
\[u_T(t,x):=u(tT,x), \ m_T(t,x):=m(tT,x), \ \theta_T(t,x):=\theta(tT,x) \quad \text{ for } (t,x) \in [0,1]\times \Omega.\]
At a formal level, we expect \cite{C13,CLLP12} that, as $T\to +\infty$,
\[ \frac1Tu_T(t,x) \sim \bar \lambda(1-t),\,\quad  \n u_T(t,x)\sim\n \bar u(x),\, \quad\, m_T \sim \bar m, \,\quad \theta_T \sim \bar\theta,\]
where $(\bar\lambda, \bar u, \bar m, \bar\theta)$ solve the so-called ``ergodic system''
\begin{equation}\label{EMFG}
        \begin{cases}
            \Bar \lambda - \nu \Delta \Bar u - \frac{\abs{\nabla \Bar u}^2}{2}= \Bar{\theta} &\quad \text{in } \ \Omega,\\
            -\nu \Delta \Bar m + \n \cdot(\Bar m  \, \nabla \Bar u)=0,  \ \bar m >0 &\quad \text{in } \ \Omega,\\
            -d \Delta \Bar \theta -\Bar \theta (k(x)-\Bar m -\Bar \theta)=0, \ \bar\theta>0 &\quad \text{in } \ \Omega,\\
            \partial_\nu \bar u = \partial_\nu \bar m = \partial_\nu \bar\theta =0 &\quad \text{on } \partial\Omega,\\
            \int_\Omega \bar u =0, \quad \fint_\Omega \bar m =1.
        \end{cases}
\end{equation}

Our main contributions are the following. We present a general theoretical result (Theorem \ref{Th:uniqueness}) that guarantees uniqueness of solution of \eqref{EMFG} in arbitrary dimension. In particular, assuming a constant $k$ (that is, the environment is homogeneous), the conditions of Theorem \ref{Th:uniqueness} imply uniqueness of steady states. Investigating the case of general $k$, Theorem \ref{Th:uniqueness_d=1} provides computable conditions on $k$ and other biologically relevant parameters guaranteeing such conditions are met in the one-dimensional case. Finally, in any dimension,  we provide in Theorem \ref{Th:attractiveniss} a condition guaranteeing global stability of the stationary state. We apply this theorem to prove that, under the assumptions of Theorem \ref{Th:uniqueness_d=1}, the stationary solution is indeed dynamically stable in the sense that solutions to the evolution system converge, in the long-time, to it. We leave the question of uniqueness of solutions to the evolution equation as an interesting open problem.

\subsection{Main model and main results}
In order to obtain a more tractable system, it is convenient to introduce the Hopf-Cole transform \[\bar\varphi := e^{\frac{\Bar u}{2 \nu}},\] which allows us to remove one of the equations in \eqref{EMFG}. Indeed, the equation for $\bar m$ can be solved explicitly, yielding $\bar m= c\bar\varphi^2$ for some constant $c\in(0,1]$. System \eqref{EMFG} can then be rewritten in the more compact form 
\begin{equation}\label{S}
    \begin{cases}
        -2\nu^2 \Delta \bar\varphi - \Bar{\theta} \bar\varphi = -\bar\lambda \bar\varphi, \ \int_\Omega \ln \varphi =0, \ \fint_\Omega c\bar\varphi^2 =1, \ \bar\varphi>0 \quad &\text{in } \Omega, \\
        -d \Delta \Bar{\theta} -\Bar{\theta} (k(x)-c\bar\varphi^2 - \Bar{\theta})=0, \quad \bar \theta>0 &\text{in } \Omega,\\
        \partial_\nu \bar \p = \partial_\nu \bar\theta = 0 &\text{on }\partial \O.
    \end{cases}
\end{equation}
In particular, when $\bar \theta$ is prescribed, $-\bar\lambda$ is also the smallest eigenvalue associated with the linear equation in \eqref{S}. From now on, we assume that
\[\fint_\Omega k > \fint_\Omega c\bar\p^2 = \fint_\Omega \bar m =1.\]
This condition guarantees \cite{BHR} that for any $m$ satisfying $\fint_\O m=1$, there exists a unique positive solution $\theta$ of 
\[-d\Delta \theta-\theta(k-m-\theta)=0.\]
From a modelling viewpoint, this condition can be interpreted as requiring that the available resources exceed the total mass of fishermen in the system.

    \subsubsection{Uniqueness results for the stationary model \eqref{EMFG}}
    Regarding the uniqueness of solutions to \eqref{EMFG}, the usual condition ensuring it is the Lasry-Lions monotonicity condition, which would write, here, as follows: for any two solutions $(\bar \lambda_i, \bar u_i, \bar m_i, \bar\theta_i)$ of \eqref{EMFG} ($i=1,2$),
    \[ \int_\O (\bar \theta_1 -\bar \theta_2)(\bar m_1 -\bar m_2) \leq 0.\] Under this condition, it is routine to check uniqueness of the solutions. Establishing such a Lasry-Lions monotonicity condition is very delicate in general. The following theorem actually shows that having enough analytical control of \emph{only one} solution of \eqref{EMFG} is sufficient to derive this inequality; the study of this analytical property is then studied depending on the specific context (\emph{e.g.} one-dimensional domains or constant $k$).
    \begin{theorem}[Sufficient condition for uniqueness of solution of \eqref{S}]\label{Th:uniqueness}
    Let $\mathcal{S}$ be the set of solutions of \eqref{EMFG}. For any $d,\nu>0$, $k \in L^\infty(\Omega)$  such that $\fint_\Omega k > 1$ and for any $(\bar \lambda_i, \bar u_i, \bar m_i, \bar\theta_i)\in \mathcal{S}$, $i=1,2$, there holds 
    \[ \int_\Omega (\bar \theta_1 -\bar \theta_2)(\bar m_1 -\bar m_2) \, dx \le  -\int_\Omega (\bar \theta_1 - \bar \theta_2)^2  \left(1- \frac{d}{4}\frac{\abs{\nabla \bar \theta_2}^2}{\bar \theta_2^3}\right).\]
    In particular, if 
    \begin{equation}\label{BaiHeLi_typeofineq_bartheta}
        d\frac{\abs{\nabla\bar\theta_2}^2}{\bar \theta_2^3} \le 4,
    \end{equation} then \eqref{EMFG} admits a unique solution.
    \end{theorem}
An important remark concerns the comparison with the coordinated framework studied in \cite{arXiv:2606.30106}. In that work, a condition similar to \eqref{BaiHeLi_typeofineq_bartheta} was shown to guarantee uniqueness of stationary solutions (see \cite[Theorem 5]{arXiv:2606.30106}). However, the arguments leading to uniqueness are different in the two settings. In the competitive scenario considered here, the structure of the running cost allows for a direct comparison between the two candidate solutions, which is the heart of the Lasry-Lions monotonicity argument. By contrast, the coordinated scenario requires an indirect approach (essentially, a linearised Lasry-Lions monotonicity).
    
    We deduce the following corollary, the proof of which is simply a combination of Theorem \ref{Th:uniqueness} and of the fact that, if $k$ is constant, $(\bar\lambda, \bar u, \bar m, \bar\theta)=(k-1, 0, 1, k-1)$ is a particular solution of \eqref{EMFG}:
    \begin{corollary}
    Assume $k>1$ is constant. \eqref{EMFG} has the unique solution $(\bar\lambda, \bar u, \bar m, \bar\theta)=(k-1, 0, 1, k-1)$.
    \end{corollary} 
    The question of whether one can find regimes where \eqref{BaiHeLi_typeofineq_bartheta} is satisfied is partially answered in the following theorem, which provides computable bounds in the one-dimensional setting. Therefore, in the remainder of this section we restrict to $\Omega=(0;\ell)$ for some $\ell>0$.
    
   \begin{theorem}[$n=1$]\label{Th:uniqueness_d=1}
        Let $k\in L^\infty((0;\ell))$ such that $\fint_0^\ell k >1$ and let $(\bar \lambda, \bar u, \bar m, \bar \theta)$ be a solution of \eqref{EMFG}. Let us also suppose that 
        \begin{equation}\label{data_condition_uniqueness}
            \frac{1+\frac{2\ell c}{\nu}\sqrt{\sup k} - c\inf k}{\fint_0^\ell k-1}e^{\frac{\ell}{\sqrt{d}}\sqrt{\sup k -\fint_0^\ell k +1}}\le 5.
        \end{equation}
        Then \eqref{BaiHeLi_typeofineq_bartheta} holds. 
    \end{theorem}

\subsubsection{Attractiveness of the stationary solution}
After studying the uniqueness of the stationary equilibrium, it is natural to investigate its dynamical stability.
    \begin{theorem}\label{Th:attractiveniss}
        Let $k\in C^1(\overline \Omega)$ satisfy $\fint_\Omega k >1$, $m_0\ge 0$ in $L^\infty(\Omega)$, $\fint_\Omega m_0=1$, and let $(\bar \lambda, \, \bar u, \, \bar m,\, \bar \theta)$ be a solution of \eqref{EMFG}. Let $\theta_0 \in L^\infty(\Omega),\, \inf_{\overline \O}\theta_0>0$ and assume that
        \begin{equation}\label{constraint}
            \forall T >0 \text{ large, a solution }(u, \, m,\, \theta) \text{ of \eqref{MFG} satisfies }\int_0^T \int_\Omega (\theta -\bar \theta) (m-\bar m) - \int_\Omega \left( \theta_0 -\bar \theta -\bar \theta \ln \left( \frac{\theta_0}{\bar \theta}\right)\right)\le 0.
        \end{equation}
        Let us consider the rescaled functions: 
        \[u_T(t,x):=u(tT,x), \ m_T(t,x):=m(tT,x), \ \theta_T(t,x):=\theta(tT,x) \quad \text{ for } (t,x) \in [0,1]\times \Omega.\]
        Then, as $T \to +\infty$, the following convergences are satisfied:
        \begin{itemize}
            \item Gradient convergence:
            \[\nabla u_T \longrightarrow \nabla \bar u \ \text{in $L^2((0;1)\times\Omega)$};\]
            \item Density convergence:
            \[m_T \longrightarrow \bar m \ \text{in $L^1((0;1)\times\Omega)$};\]
       \item Convergence of the population distribution:
            \[\theta_T \longrightarrow \bar \theta \ \text{in $L^2((0;1)\times\Omega)$ };\]
            \item Convergence of the value function: there holds
\[
   \sup_{t\in [0;1]}\left\Vert \frac{u_T(t,\cdot)}T-\bar\lambda(1-t)\right\Vert_{L^\infty( \O)}\underset{T\to\infty}\rightarrow 0\]
and 
\[
   \left\Vert u_T-\fint_\O u_T-\bar u\right\Vert_{L^2(0,1;W^{1,2}(\O))}\underset{T\to\infty}\rightarrow 0.\]                      In particular, the ergodic constant $\bar \lambda$ is the limit of the average gain for large time $\frac{\sup_\alpha J_T(0,\alpha)}{T}$.
        \end{itemize}
    \end{theorem}
    Combining Theorem \ref{Th:uniqueness} and Theorem \ref{Th:attractiveniss}, we have
    \begin{theorem}\label{Th:constraintlongtimeaverage}
    Let $k\in C^1(\overline \Omega)$ be such that $\fint_\Omega k >1$. Assume that a solution $(\bar \lambda, \, \bar u, \, \bar m,\, \bar \theta)$ of \eqref{EMFG} satisfies \eqref{BaiHeLi_typeofineq_bartheta}, so that it is the unique solution to \eqref{EMFG}. Assume also $T>0$ large, $m_0\ge 0$ in $L^\infty(\Omega)$, $\fint_\Omega m_0=1$, $\theta_0,\, \inf_{\overline \O}\theta_0>0$ in $L^\infty(\Omega)$ and consider a solution $(u, \, m,\, \theta)$ of \eqref{MFG}. Then  \eqref{constraint} holds for any $\theta_0$ and, consequently, any solution of the evolutionary system converges to the stationary solution as $T\to+\infty$. 
    \end{theorem}
    A remarkable aspect of Theorem \ref{Th:attractiveniss} and Theorem \ref{Th:constraintlongtimeaverage} is that they operate without assumptions on $m$, and do not require many constraints on initial data. We address this problem by proposing a strategy of proof inspired by \cite{CLLP13}.
    As mentioned earlier, we do not address in this paper uniqueness of solution for the evolutionary system \eqref{MFG}. We believe that such uniqueness should hold under less restrictive conditions similar to \eqref{constraint}, but that still remains an open and challenging problem.

\subsection{Bibliographical references}\label{Sec:Bibliography}
\subsubsection{Reaction Diffusion Equations} We review some basic facts about the dynamics of fishes (in the absence of fishermen); a general reference is \cite{BHR}. The fish population is modelled through the standard logistic diffusive equation
\begin{equation}\label{RDE}\tag{RDE}
    \begin{cases}
        \partial_t \theta - d \Delta \theta = \theta (g(x) -\theta) \quad &\text{in } \R_+ \times \Omega,\\
        \partial_\nu \theta =0 &\text{on } \R_+ \times \partial\Omega,\\
        \theta(0,\cdot)=\theta_0(\cdot)\not\equiv0,\ge 0 &\text{in } \Omega.
    \end{cases}
\end{equation}
 \eqref{RDE} is crucial in population dynamics \cite{S,SK97} and studied by \cite{CC, CC2} in bounded domains,  by \cite{BHR} in periodic domains, following the pioneering papers of \cite{F37, K37, S} for constant parameters. The \eqref{RDE} is also called Fisher-KPP equation; the function $\theta$ denotes a fish population, $d >0$ is a diffusivity coefficient, the reaction term is a logistic law and $g$ is the intrinsic growth rate parameter (including bird and death coefficients due to the specific species, bird and death coefficients due to the environment). The $g$ parameter represents the availability of resources:  if $g>0$ then it is a favorable zone for the fish population to survive, if $g<0$, the zone is unfavorable for the persistence of the species. The non-linearity could be interpreted as the competition between two individuals of the species to get resources. Classical results from \cite{CCH, BHR} establish that the existence of a positive nontrivial solution for \eqref{RDE}, the uniqueness and the instability of the trivial solution $\theta \equiv 0$ are directly related to the negative sign of the principal eigenvalue $\lambda_1[-d\Delta -g(x)]$ (i.e. the smallest eigenvalue) uniquely defined by the following eigenvalue problem
\[\begin{cases}
    -d\Delta \phi - g(x)\phi = \lambda \phi, \ \phi>0 \quad &\text{in } \Omega,\\
    \partial_\nu\phi = 0 &\text{on }\partial \Omega.
\end{cases}\]
Moreover, the following inequalities hold
\[ -\max_{\bar\Omega} g \le \lambda_1[-d\Delta -g(x)] \le -\fint_\Omega g\]
and so, if $\fint_\Omega g >0$, then $\lambda_1[-d\Delta-g(x)]<0$. 

\subsubsection{Mean Field Game}  Mean Field Game theory was introduced simultaneously and independently by \cite{LL0, LL1, LL2} and by \cite{CHM} to analyse models where a large number of agents interact strategically with the other players. The agents are supposed to be rational and indistinguishable, interacting through the empirical averages of quantities which depend on the state variable. In the limit $N \to +\infty$, the game may be modelled by a system of two coupled PDEs, which is named the MFG system (\cite{CDLL} for some theoretical results   and \cite{LCollegedeFrance} for a general reference). A classical (second order) MFG system is
\begin{equation} \label{MFG_CLLP}\begin{cases}
    -\partial_t u -\nu \Delta u - \frac{\abs{\nabla u}^2}{2}=F[m(t, \cdot)] \quad & \text{in } (0,T)\times \Omega,\\
    \partial_t m -\nu \Delta m +\n \cdot(m\, \nabla u) =0, \ \fint_\Omega m=1, \ m>0 & \text{in } (0,T)\times \Omega,\\
    u(T, \cdot)= 0, \ m(0,\cdot)=m_0(\cdot)>0, \ \fint_\Omega m_0=1 & \text{in } \Omega,\\
\end{cases}\end{equation}
with $F : L^2(\Omega) \longrightarrow\R$ a running cost. Although a heterogeneous final-time condition $u(T, \cdot) = G(\cdot, m(T))$ and a general Hamiltonian $H$ could  be considered, we consider a zero final-time condition for the value function $u$ and a quadratic Hamiltonian $H(p)=\frac{1}{2}p^2$ to maintain consistency with our framework.

The system \eqref{MFG} under consideration here does not enter in the classical MFG setting \eqref{MFG_CLLP}, because our running cost $\theta$ could be viewed as a functional of $m$ that is non-local both in space and time. As far as we know, most MFG systems in the literature (apart from \cite{LL0} and the perturbative MFG system in \cite{zbMATH08109759} that we discuss below) involve running costs that are either non-functional (i.e. depend directly on $m-$variable) or functionals that are non-local only in the space variable.

Let us discuss the classical results for \eqref{MFG_CLLP}. The existence of solutions for \eqref{MFG_CLLP}, whether local or non-local in space, is a straightforward consequence of a fixed point argument (see, for instance, \cite{CLLP12}). For our static \eqref{EMFG} problem, a similar approach applies to $\bar m \longmapsto \left(\bar \lambda, F[\bar m]\right)$, where $F[\bar m]=\bar\theta[\bar m]$, with minor modifications of the argument presented in the appendix of \cite{CLLP12}. In particular, one can show that the $\|F[\bar m]\|_{L^\infty(\Omega)}$ is bounded and it is independent of $\bar m$. Establishing uniqueness (or non-uniqueness) of solution for \eqref{MFG_CLLP}, however, remains a challenging problem. There are three regimes for which uniqueness of solution for \eqref{MFG_CLLP} type of systems is well-known: the so-called Lasry-Lions monotonicity condition, introduced in \cite{LL2}, which reads \begin{equation} \label{monotonicitynonlocalinspace}
    \int_\Omega (F[m_1] - F[m_2]) (m_1-m_2) \, dx \le 0 \quad \forall m_1\neq m_2.
\end{equation}
It means that a player prefers going where there are fewer other players: it could be interpreted like crowding aversion. The second setting for uniqueness is the small-time horizon case. We do not use this approach here, because we are interested in large time behavior, but let us mention it for the sake of completeness. This topic was first introduced in \cite[January 2009]{LCollegedeFrance}  and later revisited by \cite{BF}, where it is shown that the regularity of the Hamiltonian $H$ plays a crucial role in proving uniqueness for small-time horizons. In the absence of such regularity, counterexamples can be constructed, as illustrated in \cite{BF}. The last case for uniqueness is the small-data one \cite{BF, BC}. Concerning the time asymptotic for MFG systems, the first contributions were the works of Cardaliaguet, Lasry, Lions and Porretta. Under the hypothesis \eqref{monotonicitynonlocalinspace}, they showed the convergence to an ergodic MFG system in the local setting \cite{CLLP12} and in the non-local in space setting \cite{CLLP13}. For first order systems, we refer to \cite{C13}. 

\subsubsection{Harvesting Mean Field Game} The mathematical analysis of the harvesting problem within the context of Mean Field Games has seen significant recent interest and developments in \cite{MKB23, zbMATH08109759, MRB22}. In \cite{zbMATH08109759}, the authors address the following perturbative MFG system
\begin{equation}\label{perturbativeMFG}\tag{$MFG_\epsilon$}\begin{cases}
    -\partial_t u -\nu \Delta u - \frac{\abs{\nabla u}^2}{2}= \rho \star \theta_{\epsilon, m} \quad & \text{in } (0,T)\times \Omega,\\
    \partial_t m -\nu \Delta m +\n \cdot(m\, \nabla u) =0, \ \fint_\Omega m=1, \ m>0 & \text{in } (0,T)\times \Omega,\\
    \partial_t \theta_{\epsilon, m} -d \Delta \theta_{\epsilon, m} -\theta_{\epsilon, m} (k(x) -\epsilon (\rho \star m) -\theta_{\epsilon, m})=0, \ \theta_{\epsilon, m}>0, & \text{in } (0,T)\times \Omega,\\
    u(T, \cdot)= 0, \ m(0,\cdot)=m_0(\cdot)>0, \ \fint_\Omega m_0=1, \ \theta_{\epsilon, m} (0, \cdot) =\theta_0(\cdot) \ge0, \not\equiv 0 & \text{in } \Omega,\\
\end{cases}\end{equation}
where $\star$ denotes the convolution in space, $\rho$ is a smooth, even function with a compact support and such that $\rho \in C^\infty_c (\Omega)$, $\rho \ge 0$, $\int_\Omega \rho =1$.  This regime belongs to the class of small-interactions MFG and, here, the ``smallness'' is quantified by $\epsilon >0$. They  manage to derive a quantified Lasry-Lions monotonicity condition of the form
\[\int_0^T\int_\Omega (\rho \star \theta_{\epsilon,m_1} - \rho \star \theta_{\epsilon,m_2}) (m_1-m_2) \le -a \epsilon \int_0^T\int_\Omega (\theta_{\epsilon, m_1} -\theta_{\epsilon, m_2})^2, \quad a>0,\]
note that, in that case, as the non-locality is in both variables, unlike \eqref{monotonicitynonlocalinspace}, the monotonicity involves an integral in time and in space. This allowed them to derive  uniqueness and convergence to the  ergodic system. However, their regime is more restrictive than ours; indeed, in the one-dimensional case, we are able to get directly uniqueness for the ergodic case using an approach that does not rely on a perturbative argument.

In our previous paper \cite{arXiv:2606.30106}, we focus on the coordinated counterpart of the present work, namely a harvesting Mean Field Control (MFC) model. We identify $k$ as a bifurcation parameter and we show that when $k$ crosses a specific threshold, a bifurcation occurs, leading to the existence of at least two distinct ecological equilibria for the stationary MFC model. 
\section{Uniqueness results: proof of Theorem \ref{Th:uniqueness} and Theorem \ref{Th:uniqueness_d=1}} \label{Sec:Uniqueness}
    Let us start this Section with the proof of Theorem \ref{Th:uniqueness}, which provides a general, although at first glance hard to check, sufficient condition for uniqueness.
    \begin{proof}[Proof of Theorem \ref{Th:uniqueness}]
    Let $(\bar\lambda_i, \bar u_i, \bar m_i,\bar\theta_i) \in \mathcal{S}$, $i=1,2$. It is enough to observe that
        \begingroup\allowdisplaybreaks\begin{align*}
            &\int_\Omega (\bar \theta_1-\bar \theta_2)(\bar m_1-\bar m_2)\\
            &=\int_\Omega (\bar \theta_1-\bar \theta_2) \left( d\frac{\Delta \bar \theta_1}{\bar \theta_1}- \bar \theta_1 - d \frac{\Delta \bar \theta_2}{\bar \theta_2}+\bar \theta_2 \right)\\
            &=-\int_\Omega (\bar \theta_1-\bar \theta_2)^2 +d \int_\Omega \frac{\bar \theta_1-\bar \theta_2}{\bar \theta_1}\Delta \bar \theta_1 - d \int_\Omega \frac{\bar \theta_1 - \bar \theta_2}{\bar \theta_2} \Delta \bar \theta_2\\
            &= -\int_\Omega (\bar \theta_1-\bar \theta_2)^2+d \int_\Omega \left( \Delta \bar \theta_1 + \Delta \bar \theta_2 \right) - d \int_\Omega \left( \frac{\bar \theta_2}{\bar \theta_1}\Delta \bar \theta_1 + \frac{\bar \theta_1}{\bar \theta_2}\Delta \bar \theta_2 \right)\\
            &= -\int_\Omega (\bar \theta_1 - \bar \theta_2)^2 - d \int_\Omega \bar \theta_2 \bar \theta_1 \left(\frac{\Delta \bar \theta_1}{\bar \theta_1^2}+ \frac{\Delta \bar \theta_2}{\bar \theta_2^2}\right)\\
            &=-\int_\Omega (\bar \theta_1-\bar \theta_2)^2 +d \int_\Omega \nabla \bar \theta_1 \nabla \bar \theta_2 \left( \frac{1}{\bar \theta_1}+ \frac{1}{\bar \theta_2}\right) - d \int_\Omega \bar \theta_2 \bar \theta_1 \left( \frac{\abs{\nabla \bar \theta_1}^2}{\bar \theta_1^3}+ \frac{\abs{\nabla \bar \theta_2}^2}{\bar \theta_2^3}\right)\\
            &\le -\int_\Omega (\bar \theta_1 - \bar \theta_2)^2 + d\int_\Omega\frac{\abs{\bar \theta_1+\bar \theta_2}}{\bar \theta_1\bar\theta_2}\abs{\nabla \bar \theta_1} \abs{\nabla \bar \theta_2}\frac{\sqrt{\bar \theta_2}}{\sqrt{\bar \theta_2}} - d \int_\Omega \bar \theta_2 \bar \theta_1 \left( \frac{\abs{\nabla \bar \theta_1}^2}{\bar \theta_1^3}+ \frac{\abs{\nabla \bar \theta_2}^2}{\bar \theta_2^3}\right)\\
            &\le -\int_\Omega (\bar \theta_1 - \bar \theta_2)^2  + d \int_\Omega \left(\frac{\bar \theta_2\abs{\nabla \bar \theta_1}^2}{\bar \theta_1^2} + \frac{\abs{\bar \theta_1+\bar \theta_2}^2}{4}\frac{\abs{\nabla \bar \theta_2}^2}{\bar \theta_2^3}\right) - d \int_\Omega \bar \theta_2 \bar \theta_1 \left( \frac{\abs{\nabla \bar \theta_1}^2}{\bar \theta_1^3}+ \frac{\abs{\nabla \bar \theta_2}^2}{\bar \theta_2^3}\right)\\
            &= -\int_\Omega (\bar \theta_1 - \bar \theta_2)^2  \left(1- \frac{d}{4}\frac{\abs{\nabla \bar \theta_2}^2}{\bar \theta_2^3}\right).
        \end{align*}\endgroup
  The conclusion follows by assuming \eqref{BaiHeLi_typeofineq_bartheta}.
    \end{proof}
    
    Let us observe that Theorem \ref{Th:uniqueness} provides a sufficient condition for uniqueness based on \eqref{BaiHeLi_typeofineq_bartheta}, and a Bai, He \& Li-type of inequality \cite{BHL2015}. Indeed, trying to estimate $\fint_\Omega \bar \theta \le 3 \fint_\Omega (k-\bar m)$, \cite{BHL2015} derived in dimension one the following inequality
    \[ d\frac{\abs{\bar \theta'}^2}{\bar\theta^3}\le \frac{2}{3}.\]
    This fact was already noticed in \cite{zbMATH08109759}. Hence, a natural combination of these results is Theorem \ref{Th:uniqueness_d=1}.
    
    \begin{proof}[Proof of Theorem \ref{Th:uniqueness_d=1}]
         Set $\Omega =(0;\ell)$. Let us consider the Hopf-Cole transform system \eqref{S}, namely $(\bar\lambda, \bar\varphi, \bar \theta)$ solution of \eqref{S}. In particular, we have 
        \[ 0\leq d\fint_0^\ell \frac{\bar\theta'^2}{\bar\theta^2}= d\fint_0^\ell \frac{\bar\theta''}{\bar\theta}= \fint_0^\ell \bar \theta - \fint_0^\ell k +1,\]
        and so
        \[ \fint_0^\ell \bar \theta \ge \fint_0^\ell k -1.\]
        Moreover, we have
        \begin{align*}
            \ln \left( \frac{\bar\theta}{\fint_0^\ell \bar \theta}\right) &=\ln \bar \theta - \ln \left( \fint_0^\ell \bar \theta\right) \\
            &\ge - \int_0^\ell \abs{\frac{\bar\theta'}{\bar \theta}}\\
            &\ge -\ell \sqrt{\fint_0^\ell \abs{\frac{\bar\theta'}{\bar\theta}}^2} \\
            &\ge -\frac{\ell}{\sqrt{d}}\sqrt{\fint_0^\ell \bar\theta - \fint_0^\ell k +1},
        \end{align*}
        that is
        \begin{equation}\label{estimatebartheta1D}
            \bar\theta(x) \ge \left( \fint_0^\ell k -1 \right) e^{-\frac{\ell}{\sqrt{d}} \sqrt{\sup k - \fint_0^\ell k +1}} \quad \forall x\in (0,\ell). 
        \end{equation}
        Moreover, integrating the equation satisfied by $\bar \varphi$, we have
        \[-\bar\lambda \fint_0^\ell \bar\varphi = -\fint_0^\ell \bar\theta\bar\varphi \ge -\max \bar\theta \fint_0^\ell \bar\varphi,\]
        that is
        \[ \bar\lambda \le \max \bar\theta.\]
        We can also observe that 
        \[2\nu^2\fint_0^\ell\abs{\bar\varphi'}^2 =-\fint_0^\ell \left(\bar\lambda - \bar \theta\right) \bar\varphi^2 \le \frac{2\max\bar\theta}{c}  \le \frac{2\sup k}{c},\]
        and, together with the Poincar\'e-Wirtinger inequality, we obtain
        \begin{align*}
            \abs{\bar\varphi^2-\fint_0^\ell \bar\varphi^2} &\le \int_0^\ell \abs{\left( \bar\varphi^2\right)'}\\
            &= 2\ell \fint_0^\ell \abs{\bar\varphi'\bar\varphi}\\
            &\le \frac{2\ell}{\sqrt{c}}\sqrt{\fint_0^\ell \bar \varphi'^2}\\
            &\le  \frac{2\ell}{\nu}\sqrt{\sup k},
        \end{align*}
        that is 
        \begin{equation}\label{Eq:p^2_bound}
            \bar\varphi^2(x) \le \frac{1}{c}+\frac{2\ell}{\nu}\sqrt{\sup k}. 
        \end{equation}
        Let us define $v:=d(\bar\theta')^2 -\alpha \bar\theta^3$, with $\alpha \in \R^+$. There exists $\bar x \in [0,\ell]$ such that $\max_{[0,\ell]} v=v(\bar x)$. As $v$ is negative at the boundary of the domain, any positive maximum is an interior point: $\bar x \in (0,\ell)$ and so, $\bar \theta'(\bar x) \neq 0$. Let us assume by contradiction that $v(\bar x)>0$. Then
        \[ 0=v'(\bar x) =\bar\theta'(\bar x) \left[ 2d\bar\theta''(\bar x) -3 \alpha \bar\theta^2(\bar x)\right],\]
        and, since $\bar \theta'(\bar x) \neq 0$, we have
        \begin{align*} 
            \frac{3}{2}\alpha \bar \theta(\bar x)&=d \frac{\bar\theta''(\bar x)}{\bar\theta(\bar x)}\\
            &= \bar \theta(\bar x) -k(\bar x) + c\bar\varphi^2(\bar x) \\
            &\le \bar \theta(\bar x) - k(\bar x) + 1+\frac{2\ell c}{\nu}\sqrt{\sup k} -c\inf k \text{ by \eqref{Eq:p^2_bound}}.
        \end{align*}
        Thus
        \[ \left(\frac{3}2\alpha-1\right)\bar \theta(\bar x)\leq  1+\frac{2\ell c}{\nu}\sqrt{\sup k}-c\inf k. \] In particular, if 
        \[  1+\frac{2\ell c}{\nu}\sqrt{\sup k}-c\inf k\leq 5 \left(\fint_0^\ell k-1\right)e^{-\frac{\ell}{\sqrt{d}}\sqrt{\sup k -\fint_0^\ell k +1}}\] we obtain a contradiction if $\alpha>4$. Thus, for any $\alpha>4$, 
       \[ d(\bar\theta')^2 \leq \alpha \bar\theta^3\] and, letting $\alpha\to 4$, we deduce
        \[  d\frac{\bar \theta'^2(x)}{\bar \theta^3(x)}\le 4 \quad \forall x\in(0,\ell).\]
        Theorem \ref{Th:uniqueness} now applies.
    \end{proof}
    Note that, instead of using \eqref{estimatebartheta1D}, one could alternatively observe, by the maximum principle, that
    \[\bar\theta\geq \inf\bar\theta \geq \inf k-c\|\bar\varphi\|_{L^\infty((0;\ell))}^2 .\]
    In this case, the estimate \eqref{BaiHeLi_typeofineq_bartheta} would follow for every $\alpha\geq \frac23$ provided that
    \[\inf k\geq \frac{1}{c}+\frac{2\ell}{\nu}\sqrt{\sup k}.\]
    Consequently, one would directly obtain
    \[d\frac{\bar\theta'^2(x)}{\bar\theta^3(x)}\leq \frac23,\qquad x\in(0,\ell),\]
    which is enough to apply Theorem \ref{Th:uniqueness}, since $\frac23\leq 4$. Hence, in this alternative approach, assumption \eqref{data_condition_uniqueness} is not required. On the other hand, estimate \eqref{estimatebartheta1D} allows us to reach the same conclusion without imposing the additional condition
    \[\inf k\geq \frac{1}{c}+\frac{2\ell}{\nu}\sqrt{\sup k}.\]
    
    Concerning the dimensional restriction, our main inspiration for the previous result comes from the \cite{BHL2015} estimate. Nevertheless, by adapting the technique developed in \cite{CLPW}, we expect that the result can be extended to arbitrary dimensions.
    
\section{Attractiveness results: proof of Theorem \ref{Th:attractiveniss} and Theorem \ref{Th:constraintlongtimeaverage}}\label{Sec:Attractiveness}
The aim of this Section is to prove the attractiveness of the unique stationary solution under the ``monotonicity'' assumption \eqref{constraint}. Let us observe that if $(u,m,\theta)$ is a solution of \eqref{MFG}, then the rescaled functions $(u_T,m_T, \theta_T)$ solve
    \begin{equation}\label{MFGT}
        \begin{cases}
            -\frac{1}{T}\partial_t u_T -\nu \Delta u_T -\frac{\abs{\nabla u_T}^2}{2}= \theta_T &\text{in } \ (0;1) \times \Omega,\\
            \frac{1}{T}\partial_t m_T -\nu \Delta m_T +\n \cdot(m_T \nabla u_T)=0, \ \fint_\Omega m_T =1 \quad &\text{in } \ (0;1)\times \Omega,\\
            \frac{1}{T}\partial_t \theta_T -d \Delta \theta_T - \theta_T (k-m_T-\theta_T)=0, \ \theta_T >0 \quad & \text{in } \ (0;1) \times \Omega, \\
            \partial_\nu u_T= \partial_\nu m_T= \partial_\nu \theta_T=0 & \text{on }  (0;1) \times \partial\Omega, \\
            u_T(1,x)=0,  \ m_T(0,x)=m_0(x), \  \theta_T(0,x)=\theta_0(x),\, \inf_{\overline\O}\theta_0>0 \quad &\text{in } \ \Omega.
        \end{cases}
    \end{equation}
    The overall strategy follows the one developed in \cite{CLLP12}. In this earlier paper, the starting point is the energy identity 
    \begin{equation}\label{Eq:keyquantity}
        -\int_0^1\int_\Omega \left(\theta_T-\bar\theta\right) \, (m_T-\bar m) + \int_0^1 \int_\Omega \frac{m_T+\bar m}{2} \abs{\nabla u_T - \nabla \bar u}^2 = \left[ \frac{1}{T} \int_\Omega (u_T-\bar u) \, (m_T- \bar m) \right]_0^1,
    \end{equation}
    from which the Lasry-Lions monotonicity condition provides control of 
    \[\Vert \nabla u_T -\nabla \bar u\Vert_{L^2((0;1)\times\Omega)}\]
    in terms of the quantity on the right-hand side of \eqref{Eq:keyquantity}. In the framework of \cite{CLLP12}, the local structure and the regularity of the running cost, which admits a time-independent primitive, yield a uniform bound on $\Vert \nabla u_T(0)\Vert_{L^2(\Omega)}$. Consequently, the right-hand side of \eqref{Eq:keyquantity} converges to zero as $T \to +\infty$. This implies the strong convergence of $\nabla u_T \to \nabla \bar u$, from which the convergence of $m_T$ and subsequently of $\frac{u_T}{T}$ follows.\\
    
    In our setting, this argument cannot be applied directly. First, the boundedness of the stationary density $\bar m$ is not available a priori and must be established by regularity estimates (Lemma \ref{Lem:regularitytheta}). Second, the running cost is non-local and it no longer enjoys the regularity presented in \cite{CLLP12}, so a different argument is required (Lemma \ref{Lem:boundnablauT(0)}). Finally, the non-local dependence of the running cost allows us to introduce a correction to the classical Lasry-Lions monotonicity condition, namely the condition \eqref{constraint}, leading to the same conclusion under a weaker assumption (Proposition \ref{Prop}).\\
    
    Let us start with some preliminary results.
    \begin{lemma}\label{Lem:boundnablauT(0)}
        Let $(u_T, m_T, \theta_T)$ be a solution of \eqref{MFGT} and assume that $\theta_T \in L^\infty((0;1)\times\Omega)$. In particular, the function $u_T$ solves
        \[ \begin{cases}
            -\frac{1}{T}\partial_t u_T - \nu \Delta u_T - \frac{\abs{\nabla u_T}^2}{2}=\theta_T &\quad \text{in } \ (0;1)\times \Omega,\\
            u_T(1,\cdot)=0 &\quad \text{in } \ \Omega.
        \end{cases}\]
        Then for any $t\in [0,1]$ \[\mathrm{osc}(u_T)(t,\cdot):=\sup_{\overline \O}u_T-\inf_{\overline\O}u_T\leq C(\nu,\left\Vert \theta_T\right\Vert_{L^\infty((0;1)\times \O)}, \delta) .\]
    \end{lemma}
    
    \begin{proof}[Proof of Lemma \ref{Lem:boundnablauT(0)}]
  We use once more the Hopf-Cole transform. Indeed, the function $\phi_T:=e^{\frac{u_T(1-\cdot,\cdot)}{2\nu}}$ solves
  \[ \frac1T \partial_t\phi_T-\nu\Delta \phi_T=\theta_T\frac{\phi_T}{2\nu},\,\quad \phi_T(0,\cdot)=1.\] 
  Hence, for any $\delta>0$, $\varphi_T (t,x):=\phi_T (Tt,x)$ is a solution of the extended equation
    \[ \partial_t\varphi_T-\nu\Delta \varphi_T=\theta_T\frac{\varphi_T}{2\nu}{\bf 1}_{t\geq 0} \hbox{ in } (-1,T),\,\quad \varphi_T(-1,\cdot)=1.\] 

  By the Harnack inequality \cite{zbMATH05045580} (see also \cite[Theorem 1.2.6]{zbMATH07668634}), for some constant $c=c(\nu,\left\Vert \theta_T\right\Vert_{L^\infty((0;1)\times \O)})$, and for all $t\geq 0$:
  \[ \sup_{x\in \overline \O} \varphi_T(t,x)\leq c\inf_{x\in \overline \O}\varphi_T(t,x).\] Coming back to $\phi_T$ and passing to the logarithm gives the conclusion.
    \end{proof}
    
    \begin{lemma}\label{Lem:regularitytheta}
    Let $\Omega \subset \R^n$ be a smooth bounded domain, and let $k\in L^\infty(\Omega)$ satisfy $\fint_\Omega k>1$. Let $(u_T,m_T,\theta_T)$ be a solution of \eqref{MFGT}, and let $(\bar\lambda,\bar u,\bar m,\bar\theta)$ be a solution of \eqref{EMFG}. Then there exists a constant $C>0$, independent of $T$, such that
    \[0\le \theta_T(t,x)\le C\quad\text{for all }(t,x)\in (0,1)\times\Omega,\]
    and
    \[0< \bar\theta(x)\le C \quad\text{for all }x\in\Omega.\]
    Moreover, we have the following regularity for the stationary solution 
    \[\bar u, \bar\theta \in C^{1,\alpha}(\overline\Omega), \quad \bar m \in C^{2,\alpha}(\overline\Omega) \ \text{ for some } \alpha\in (0;1).\]
    If, in addition, $k\in C^{1,\alpha}(\overline\Omega)$, then standard bootstrap arguments imply $\bar u, \bar\theta \in C^{2,\alpha}(\overline\Omega)$.
    \end{lemma}
    \begin{proof}
        Since $\fint_\Omega k >1$, we have $\lambda_1\left[ -d\Delta - (k-\bar m)\right]<0$. Since $\bar m > 0$, it follows immediately that $\lambda_1\left[-d\Delta-k\right]<0$. Let $\bar \theta_k$ denote the unique positive solution of 
        \[\begin{cases}
            -d\Delta \bar \theta_k- \bar \theta_k (k(x)-\bar\theta_k)=0, \ \bar\theta_k>0 &\quad \text{in } \ \Omega,\\
            \partial_\nu\bar \theta_k =0 &\quad \text{on } \ \partial\Omega.
        \end{cases}\]
        
        Let us consider $\bar \theta_k \ge \max_\Omega \theta(0,\cdot)>0$. We observe that $\sup k$ is a supersolution of this problem. Moreover, since $m_T, \, \bar m >0$, the same function $\sup k$ is also a supersolution for both $\bar \theta$ and $\theta_T$. By the comparison principle, we obtain
        \[ 0\le \theta_T(t,x)  \le \max \left\{ \|\theta_0\|_{L^\infty(\Omega)},\sup k\right\} \quad \forall(t,x) \in (0;1)\times\Omega, \ \forall T>0\]
        and \[ 0<\bar \theta(x) \le \sup k \quad \forall x \in \Omega.\]
        Let us consider the Hopf-Cole transform $\bar\p=e^{\frac{\bar u}{2\nu}}$, so that $\bar m=c\bar\p^2$ for some constant $c\in(0,1]$ and $(\bar\lambda,\bar\p,\bar\theta)$ solves
        \[\begin{cases}
            -2\nu^2\Delta \bar\p= (\bar\theta-\bar\lambda)\bar\p, \ \int_\Omega \ln \varphi =0,\ \bar\p>0, \ \fint_\Omega c\bar\p^2=1 \quad &\text{in } \Omega,\\
            -d\Delta\bar\theta = \bar\theta(k- c\bar\p^2-\bar\theta),\bar\theta>0  &\text{in } \Omega,\\
            \partial_\nu \bar\p =\partial_\nu \bar\theta=0  &\text{on } \partial\Omega.\\
        \end{cases}\]
        Integrating the first equation yields
        \[\bar\lambda \int_\Omega \bar\p = \int_\Omega \bar\theta \bar\p,\]
        hence $\bar\lambda  \le \sup k$ and so, $\Vert \bar\theta-\bar\lambda\Vert_{L^\infty} \le 2\sup k$.\\
The regularity of the solutions follows from standard elliptic estimates.

    \end{proof}
    
    Let us recall, for the sake of completeness, the statement of a preliminary identity from \cite{CLLP12}.
    \begin{lemma}[\cite{CLLP12}, Lemma 2.2]\label{Lem:keyidentity}
        For any $0\le t_1<t_2\le 1$ we have
        \begin{equation*}
            \left[ \int_\Omega (u_T-\bar u) (m_T-\bar m) \, dx \right]^{t_2}_{t_1} + T \int^{t_2}_{t_1} \int_\Omega \left[\frac{m_T+\bar m}{2}\abs{\nabla u_T -\nabla \bar u}^2 -\left(\theta_T-\bar\theta\right) \left( m_T-\bar m\right) \right] \, dxdt=0.
        \end{equation*}
    \end{lemma}
    
    \begin{proposition}\label{Prop}
    Let $\Omega \subset \R^n$ be a smooth bounded domain, and let $k\in L^\infty(\Omega)$ satisfy $\fint_\Omega k>1$. Let $(u_T,m_T,\theta_T)$ be a solution of \eqref{MFGT}, and let $(\bar\lambda,\bar u,\bar m,\bar\theta)$ be a solution of \eqref{EMFG}. Let $\theta_0 \in L^\infty(\Omega),\, \inf_{\overline \O}\theta_0>0$ and assume that
    \begin{equation}\label{assumption}
        \int_0^1\int_\Omega (\theta_T-\bar \theta) (m_T-\bar m) - \frac{1}{T}\int_\Omega \left[ \theta_0-\bar \theta - \bar \theta \ln \left(\frac{\theta_0}{\bar \theta}\right)\right] \le 0.
    \end{equation}
    Then 
    \begin{equation}\label{L2estimatenablauT-nablabaru}
        \Vert \nabla u_T -\nabla \bar u\Vert_{L^2((0;1)\times\Omega)}\sim O\left(\frac1{\sqrt{T}}\right)
    \end{equation}
    and
    \begin{equation}\label{Eq.mT(nabla uT nabla u)2} 
        \lim_{T \to +\infty} \int_0^1\int_\Omega m_T \abs{\nabla u_T - \nabla \bar u}^2=0.
    \end{equation}
    \end{proposition}
    \begin{proof}
    Proceeding as in \cite[Lemma 2.4]{CLLP12} and using Lemma \ref{Lem:keyidentity}, we obtain
    \begin{align*}
        \int_0^1\int_\Omega \left(\theta_T-\bar\theta\right) \, (m_T-\bar m) &- \int_0^1 \int_\Omega \frac{m_T+\bar m}{2} \abs{\nabla u_T - \nabla \bar u}^2\\
        &= -\left[ \frac{1}{T} \int_\Omega (u_T-\bar u) \, (m_T- \bar m) \right]_0^1 \\
        &= \frac{1}{T} \left[\int_\Omega \bar u (m_T(1,\cdot)-\bar m) - \int_\Omega (u_T(0,\cdot)-\bar u) (m_0-\bar m) \right].
    \end{align*}
    Since $\bar u \in L^\infty(\Omega)$ by Lemma \ref{Lem:regularitytheta}, and $\fint_\O m_T(t,\cdot)=\fint_\O \bar m=1$,
    \begin{equation*}
        \abs{\int_\Omega  \bar u (m_T(1,\cdot)-\bar m)} \le 2 \Vert \bar u\Vert_{L^\infty(\Omega)}.
    \end{equation*}

By Lemma \ref{Lem:boundnablauT(0)}, there exists some $C_2>0$  that does not depend on $T$ such that:
    \begin{equation*}\begin{aligned}
        \abs{\int_\Omega (u_T(0,\cdot)-\bar u) (m_0-\bar m)} 
        &\le C\left(\mathrm{osc}(u_T(0,\cdot))+\mathrm{osc}(\bar u)\right)\\&\leq C_2.
    \end{aligned}\end{equation*}

    From Lemmas \ref{Lem:boundnablauT(0)}--\ref{Lem:regularitytheta}, we obtain
    \begin{equation}
        \abs{\left[ \frac{1}{T} \int_\Omega (u_T-\bar u) \, (m_T- \bar m) \right]_0^1} \le \frac{C_1}{T}.
    \end{equation}
    Adding the non-negative correction term
    \[ \frac{1}{T} \int_\Omega \left[ \theta_0-\bar \theta -\bar \theta \ln \left(\frac{\theta_0}{\bar \theta}\right) \right]\]
    and using assumption \eqref{assumption}, we obtain
    \[\int_0^1\int_\Omega \frac{m_T+\bar m}{2} \abs{\nabla u_T -\nabla \bar u}^2\sim O(T^{-1}),\]
    and therefore, since $\inf \bar m>0$ ,
    \[ \Vert \nabla u_T -\nabla \bar u\Vert_{L^2((0;1)\times\Omega)}\sim O\left(\frac1{\sqrt{T}}\right).\]
Furthermore,
    \[ \int_0^1 \int_\Omega m_T \abs{\nabla u_T -\nabla\bar u}^2 \underset{T \to +\infty}{\longrightarrow} 0.\]
    \end{proof}

    We can move on to the convergence of densities $m_T$; to do so, we rely on an entropy argument inspired by the JKO scheme \cite{zbMATH01150904}.
    \begin{proposition}\label{Pr:ConvergenceFishermen}
    There holds    \[ m_T\underset{t\to \infty}\rightarrow \bar m\text{ in }L^1((0;1)\times \O).\]
    \end{proposition}
    \begin{proof}[Proof of Proposition  \ref{Pr:ConvergenceFishermen}]
    We introduce the relative entropy functional
    \[ S_T:[0;1]\ni  t\mapsto \int_\O \left(m_T\ln\left(\frac{m_T}{\bar m}\right)-m_T+\bar m\right).\]
    By direct computations, we obtain 
    \begin{align*}
   \frac1T \frac{dS_T}{dt}&=\frac1T\int_\O \partial_t m_T\ln\left(\frac{m_T}{\bar m} \right)
   \\&=-\nu \int_\O m_T\left|\n\ln\left(\frac{m_T}{\bar m}\right)\right|^2
   \\&+\int_\O m_T\left(\n u_T-\n \bar u\right)\cdot\n \ln\left(\frac{m_T}{\bar m}\right)
   \\&\leq -\frac{\nu}2\int_\O m_T\left|\n \ln\left(\frac{m_T}{\bar m}\right)\right|^2+\frac1{2\nu}\int_\O m_T\left|\n u_T-\n \bar u\right|^2.
    \end{align*}
    We thus deduce that 
    \[ \int_0^1\int_\O m_T\left|\n \ln\left(\frac{m_T}{\bar m}\right)\right|^2\underset{T\to \infty}\rightarrow 0. \] As $\inf \bar m>0$, there exists a constant $ C$ such that for any $\varphi\in W^{1,2}(\O)$ there holds
    \[ \int_\O \left(\varphi-\fint_\O \varphi \bar m\right)^2\bar m\leq C \int_\O |\n \varphi|^2\bar m.\] In particular, applying this inequality to $\varphi:=\frac{\sqrt{m_T}}{\sqrt{\bar m}}$, we obtain
    \[ \int_0^1\int_\O \left(\varphi-\fint_\O \varphi \bar m\right)^2\bar m\leq C\int_0^1\int_\O |\n \varphi|^2\bar m.\] Given that 
    \[ \int_0^1\int_\O |\n \varphi|^2\bar m=\frac14 \int_0^1\int_\O \left|\n \ln\left(\frac{m_T}{\bar m}\right)\right|^2m_T\] we deduce that 
    \[  \int_0^1\int_\O \left(\varphi-\fint_\O \varphi \bar m\right)^2\bar m\leq C \int_0^T\int_\O |\n \varphi|^2\bar m\underset{T\to \infty}\rightarrow 0.\]
    Now, notice that 
    \begin{align*}
      \Vert m_T-\bar m\Vert_{L^1((0;1)\times \O)}&=\int_0^1\int_\O |\p^2-1|\bar m
      \\&=\int_0^1\int_\O |\p-1|\cdot|\p+1|\bar m
      \\&\leq C \int_0^1 \sqrt{\int_\O \bar m|\p-1|^2}\sqrt{\int_\O \bar m|\p+1|^2}
      \\&=C \int_0^1 \sqrt{2-2\int_\O \bar m\p}\sqrt{2+2\int_\O \bar m\p}
      \\&=C\int_0^1 \sqrt{\left(\int_\O \bar m\p^2-\left(\int_\O \bar m\p\right)^2\right)}
      \\&\leq C\int_0^1\sqrt{\fint_\O \bar m\left(\p-\fint_\O \bar m\p\right)^2}
    \end{align*}
    and so 
\[
    \Vert m_T-\bar m\Vert_{L^1((0;1)\times \O)}^2\underset{T\to \infty}\rightarrow 0.\]
    \end{proof}
 
     We now prove the following:
     \begin{proposition}\label{Pr:ConvergenceTheta}
     There holds
     \[ \Vert \theta_T-\bar \theta\Vert_{L^2((0;1)\times \O)}\underset{T\to \infty}\rightarrow 0.\]
     \end{proposition}
    \begin{proof}[Proof of Proposition \ref{Pr:ConvergenceTheta}]
    We introduce
    \[ z_T:=\frac{\theta_T}{\bar\theta} \] and, as $\inf_{\overline \O}\bar\theta>0,\, \sup_{\overline\O}\bar\theta<\infty$, it suffices to prove that $z_T\underset{T\to \infty}\rightarrow 1$ in $L^2((0;1)\times \O)$.
    Observe that $z_T$ solves
    \begin{equation}\label{Eq:zt}\frac1T\partial_tz_T-\frac{d}{\bar\theta^2}\n\cdot(\bar\theta^2\n z_T)=z_T\left(\bar\theta(1-z_T)-(m_T-\bar m)\right).\end{equation}Introduce the relative entropy 
    \[ S(t):=\int_\O \bar\theta^2(z_T-1-\ln(z_T)),\] so that 
    \[ \frac1T\frac{dS}{dt}=\int_\O \frac1T\partial_tz_T \left(1-\frac1{z_T}\right)\bar\theta^2.\] Considering this entropy is inspired by \cite{arXiv:1308.6471}. Using \eqref{Eq:zt} yields
    \begin{align*}
    \frac1T\frac{dS}{dt}&=-d\int_\O \bar\theta^2|\n\ln( z_T)|^2-\int_\O \bar\theta^3(1-z_T)^2+\int_\O \bar\theta^2(1-z_T)(m_T-\bar m).
    \end{align*}
    This gives in particular 
    \[ \int_0^1\int_\O \bar \theta^3(1-z_T)^2\leq \frac{S(0)-S(1)}T+C \Vert m_T-\bar m\Vert_{L^1((0;1)\times \O)}\] for some constant $C$. The conclusion follows.
    \end{proof}We now prove the following result:
    \begin{proposition}\label{Pr:ConvergenceU}
    There holds
    \begin{equation}\label{Eq:CVG1}
   \sup_{t\in [0;1]}\left\Vert \frac{u_T(t,\cdot)}T-\bar\lambda(1-t)\right\Vert_{L^\infty( \O)}\leq C\left(\frac1T+\Vert \theta_T-\bar\theta\Vert_{L^2((0;1)\times \O)}\right)
    \end{equation}
and 
    \begin{equation}\label{Eq:CVG2}
    \left\Vert u_T-\fint_\O u_T-\bar u\right\Vert_{L^2(0,1;W^{1,2}(\O))}\leq C\left(\frac1{\sqrt{T}}+\Vert \theta_T-\bar\theta\Vert_{L^2((0;1)\times \O)}\right).
    \end{equation}    
    \end{proposition}
    \begin{proof}[Proof of Proposition \ref{Pr:ConvergenceU}]
    We begin by proving \eqref{Eq:CVG1}. To this end, introduce
    \[ q_T:=\frac{u_T-\bar u}T-\bar\lambda(1-t),\] which satisfies
    \[
    \begin{cases}
    -\partial_t q_T-T\left(\nu\Delta q_T+\langle \n \bar u,\n q_T\rangle\right)-\frac{T^2}2|\n q_T|^2=\theta_T-\bar \theta,
\\ q_T(1,\cdot)=-\frac{\bar u}T\end{cases}\] endowed with Neumann boundary conditions. Multiplying by $\bar m$ and integrating by parts, we deduce
\[ -\frac{d}{dt}\int_\O q_T\bar m=\int_\O \bar m(\theta_T-\bar \theta)+\frac12\int_\O |\n u_T-\n \bar u|^2\bar m.
\] Integrating in time, and using the fact that $\int_\O q_T(1,\cdot)\bar m=O\left(\frac1T\right)$ we deduce from Propositions  \ref{Prop}--\ref{Pr:ConvergenceTheta} that 
\[ \sup_{t\in [0;1]}\left| \int_\O q_T \bar m\right|\underset{T\to \infty}\rightarrow 0.\] Now, recall that from Lemma \ref{Lem:boundnablauT(0)} $\mathrm{osc}(u_T)$ is uniformly bounded in time. Thus
\[ \sup_{t\in [0;1]}\left\Vert q_T(t,\cdot)\right\Vert_{L^\infty(\O)}\leq C\mathrm{osc}(q_T)+\left| \int_\O q_T\bar m\right|\underset{T\to\infty}\rightarrow 0,\] and this establishes the convergence.

The estimate \eqref{Eq:CVG2} follows directly from \eqref{Eq:CVG1} and Proposition \ref{Prop}.
    \end{proof}
    We are now ready to prove Theorem \ref{Th:attractiveniss}.
    \begin{proof}[Proof of Theorem \ref{Th:attractiveniss}]
The theorem follows from Propositions \ref{Prop}--\ref{Pr:ConvergenceFishermen}--\ref{Pr:ConvergenceTheta}--\ref{Pr:ConvergenceU}.

    \end{proof}

    We conclude this section by presenting the proof of Theorem \ref{Th:constraintlongtimeaverage}.
    \begin{proof}[Proof of Theorem \ref{Th:constraintlongtimeaverage}]
        Under the assumptions of Theorem \ref{Th:constraintlongtimeaverage}, we know by Theorem \ref{Th:uniqueness} that there exists a unique solution $(\bar \lambda, \bar u, \bar m, \bar \theta)$ to \eqref{EMFG}.\\
         
        We begin by expressing $m$ and $\bar m$ in terms of $\theta$ and $\bar \theta$:
        \[m(t,x)=-\frac{\partial_t \theta -d \Delta\theta}{\theta}- \theta +k, \quad \bar m(x)= \frac{d\Delta \bar \theta}{\bar \theta}-\bar \theta +k.\]
        After standard manipulations, we obtain the identity
        \begin{align*}
            \int_0^T\int_\Omega (\theta-\bar \theta) (m-\bar m)&= -d \int_0^T\int_\Omega \left( \frac{\abs{\nabla \theta}^2}{\theta^2}\bar \theta + \frac{\abs{\nabla \bar \theta}^2}{\bar \theta^2}\theta \right) +d \int_0^T\int_\Omega \nabla \bar \theta \cdot \nabla \theta \left( \frac{1}{\theta } + \frac{1}{\bar \theta}\right) +\\
            &\quad - \int_0^T\int_\Omega \frac{\theta - \bar \theta}{\theta} \partial_t  \theta -\int_0^T\int_\Omega (\theta-\bar \theta)^2\\
        \end{align*}
        and, by the Young inequality and Theorem \ref{Th:uniqueness},
        \begin{align*}
            &\int_0^T\int_\Omega (\theta-\bar \theta) (m-\bar m)\\
            &\le -\int_0^T\int_\Omega d \left( \frac{\abs{\nabla \theta}^2}{\theta^2}\bar \theta + \frac{\abs{\nabla \bar \theta}^2}{\bar \theta^2}\theta \right) +d \int_0^T\int_\Omega \frac{\abs{\nabla \theta}^2}{\theta^2}\bar \theta + d \int_0^T\int_\Omega \frac{\abs{\nabla \bar \theta}^2}{4\bar \theta^2} \frac{\abs{\theta + \bar \theta}^2}{\bar \theta}  -\int_0^T\int_\Omega (\theta-\bar \theta)^2 \\
            &\qquad - \int_\Omega \left[ \theta(T) -\theta_0 -\bar \theta \ln\left(\frac{\theta(T)}{\theta_0}\right)\right]\\
            &= -\int_0^T\int_\Omega d \frac{\abs{\nabla \bar \theta}^2}{4\bar \theta^3} \left( 4\theta\bar \theta - (\theta + \bar \theta)^2\right) - \int_0^T\int_\Omega (\theta-\bar \theta)^2 - \int_\Omega \left[ \theta(T)-\theta_0 - \bar \theta \ln\left(\frac{\theta(T)}{\theta_0}\right)\right]\\
            &\le -(1-\alpha) \int_0^T\int_\Omega (\theta-\bar\theta)^2 - \int_\Omega \left[ \theta(T)-\theta_0 - \bar\theta \ln\left( \frac{\theta(T)}{\theta_0}\right)\right],
        \end{align*}
        with $\alpha \le 1$.\\
        
        Subtracting $\int_\Omega \left( \theta_0-\bar \theta - \bar \theta \ln{\left(\frac{\theta_0}{\bar \theta}\right)}\right)$ from both sides and reorganising the terms, we obtain
        \[\int_0^T\int_\Omega (\theta-\bar \theta)(m-\bar m) - \int_\Omega \left( \theta_0-\bar \theta - \bar \theta \ln{\left(\frac{\theta_0}{\bar \theta}\right)}\right) \le -(1-\alpha) \int_0^T \int_\Omega (\theta-\bar\theta)^2 - \int_\Omega \left[ \theta(T) - \bar \theta - \bar \theta \ln \left( \frac{\theta(T)}{\bar \theta}\right)\right].\]
        The last term is non-positive since the function $f(x):= x-\bar \theta -\bar \theta \ln{\left( \frac{x}{\bar \theta}\right)}$ is convex and reaches its minimum at $x=\bar \theta$.
        Finally, passing to the rescaled problem \eqref{MFGT}, let $\theta_T(t,x)=\theta(tT,x)$ and $m_T(t,x)=m(tT,x)$, we obtain
        \[ \begin{aligned}
            &\int_0^1 \int_\Omega (\theta_T -\bar \theta) (m_T-\bar m) - \frac{1}{T} \int_\Omega \left[ \theta_0 -\bar \theta - \bar \theta \ln \left(\frac{\theta_0}{\bar \theta}\right)\right] \\
            &\le -\left(1-\alpha\right) \int_0^1\int_\Omega (\theta_T -\bar \theta)^2 - \frac{1}{T} \int_\Omega \left[ \theta_T(1)-\bar\theta-\bar\theta \ln\left( \frac{\theta_T(1)}{\bar\theta}\right)\right] \le 0.
        \end{aligned}\]
        This shows that condition \eqref{constraint} holds. Hence, by Theorem \ref{Th:attractiveniss}, the unique solution of \eqref{EMFG} is attractive for the solution of the evolutionary \eqref{MFG} system.
    \end{proof}

\bibliographystyle{abbrv}
\bibliography{biblio}

\end{document}